\documentclass{amsart}

\usepackage{url}
\usepackage{amsfonts}
\usepackage{booktabs}

\numberwithin{equation}{section}

\begin{document}

\title[An octic diophantine equation and related elliptic curves]{An octic diophantine equation and related families of elliptic curves}

\author{Ajai Choudhry \and Arman Shamsi Zargar}
\address{13/4 A Clay Square, Lucknow - 226001, INDIA}
\email{ajaic203@yahoo.com}
\address{Department of Mathematics and Applications, University of Mohaghegh Ardabili, Ardabil, IRAN}
\email{zargar@uma.ac.ir}


\begin{abstract}
We obtain two parametric solutions of the diophantine equation $\phi(x_1, x_2, x_3)=\phi(y_1, y_2, y_3)$ where
$\phi(x_1, x_2, x_3)$ is the octic form defined by $\phi(x_1, x_2, x_3)=x_1^8+ x_2^8 + x_3^8 - 2x_1^4x_2^4 - 2x_1^4x_3^4 - 2x_2^4x_3^4$. These parametric solutions yield infinitely many examples of two equiareal triangles whose sides are perfect squares of integers. Further, each of the two parametric solutions leads to a family of elliptic curves of rank~$5$ over $\mathbb{Q}(t)$. We study one of the two  families in some detail and determine a set of five free generators for the family.
\end{abstract}

\maketitle

\section{Introduction}\label{Sec1}

It is well-known that it is very difficult to obtain parametric solutions of higher degree diophantine equations. In fact, the existing literature does not seem to contain any examples of parametric solutions of diophantine equations of the type,
\begin{equation*}
f(x_1, x_2, \ldots, x_n)=f(y_1, y_2, \ldots, y_n), \label{deqngen}
\end{equation*}
where $f(x_1, x_2, \ldots, x_n)$ is a form of degree $8$ in $n$ variables and $n \leq 6$.

Let $\phi(x_1, x_2, x_3)$ be the symmetric octic form in $3$ variables defined by
\begin{equation*}
\phi(x_1, x_2, x_3)=x_1^8+ x_2^8 + x_3^8 - 2x_1^4x_2^4 - 2x_1^4x_3^4 - 2x_2^4x_3^4.
\end{equation*}
In this paper we obtain two parametric solutions of the diophantine equation,
\begin{equation}
\phi(x_1, x_2, x_3)=\phi(y_1, y_2, y_3), \label{octiceqn}
\end{equation}
and show how infinitely many parametric solutions of Eq.~\eqref{octiceqn} may be obtained.

The diophantine equation~\eqref{octiceqn} arises naturally if we consider the problem of finding two equiareal triangles whose sides are perfect squares of integers. The area $A$ of a triangle, whose sides have lengths $a, b $ and $c$, is given by the well-known Heron's formula,
\begin{equation}
A=\sqrt{s(s-a)(s-b)(s-c)}, \label{Heron}
\end{equation}
where $s$ is the semi-perimeter of the triangle, that is, $s=(a+b+c)/2$.  We may write formula~\eqref{Heron} equivalently as follows:
\begin{equation}
A=\sqrt{-(a^4 + b^4+ c^4- 2a^2b^2  - 2a^2c^2 - 2b^2c^2)}/4. \label{Heron2}
\end{equation}

If we consider two triangles whose sides are  $x_i^2$ and $y_i^2$, $i=1, 2, 3$, and we use formula~\eqref{Heron2}  for the area of a triangle, the condition that the two triangles have equal areas is expressed precisely by the octic diophantine equation~\eqref{octiceqn}. It follows that parametric solutions of Eq.~\eqref{octiceqn} immediately yield examples of two equiareal triangles whose sides are perfect squares of integers. In this context, we note that  Choudhry and Zargar~\cite{Ch-Sh} have already obtained numerical examples of two equiareal and equiperimeter triangles whose sides are perfect squares, but a parametric solution of the problem was not found.

Further, integer solutions of the octic equation~\eqref{octiceqn} immediately yield six rational points (not necessarily distinct) on an  elliptic curve defined by the equation,
\begin{equation}
Y^2=X^3+DX, \label{ecDgen}
\end{equation}
where $D$ is a certain rational number. In this context it is pertinent to note that several authors~\cite{IN, Ma, Nag} have studied elliptic curves of type~\eqref{ecDgen}. In fact, Izadi and Nabardi~\cite{IN} have obtained a family of such elliptic curves of rank~$\geq 5$ using numerical solutions of Eq.~\eqref{octiceqn}.

We note that if we take $D=\phi(x_1, x_2, x_3)/4$ in Eq.~\eqref{ecDgen}, that is, we consider the elliptic curve,
\begin{equation}
Y^2=X^3+\frac{\phi(x_1, x_2, x_3)}{4}X, \label{cubicec}
\end{equation}
defined over the function field  $\mathbb{Q}(x_1, x_2, x_3)$, we immediately get three points $P_1$, $P_2$, $P_3$, on the curve~\eqref{cubicec}, where $P_i=(X_i, Y_i)$, $i=1, 2, 3$, are given  by
\begin{equation*}
\begin{aligned}
P_1&=\big(x_1^2x_2^2, x_1x_2(x_1^4 + x_2^4 - x_3^4)/2\big),\\
P_2&=\big(x_1^2x_3^2, x_1x_3(x_1^4 - x_2^4 + x_3^4)/2\big),\\
P_3&=\big(x_2^2x_3^2, x_2x_3(x_1^4 - x_2^4 - x_3^4)/2\big).
\end{aligned}
\label{pointsP}
\end{equation*}

When the diophantine equation~\eqref{octiceqn} is satisfied, we readily get three more points, with abscissae $y_1^2y_2^2$, $y_1^2y_3^2$, $y_2^2y_3^2$, on the curve~\eqref{cubicec}. On using a parametric solution of the octic equation~\eqref{octiceqn}, we get a parametric family of elliptic curves of type~\eqref{ecDgen} whose generic rank is $5$. We will study this family of elliptic curves in greater detail in Section~\ref{Sec3} of this paper.

\section{The octic diophantine equation $\phi(x_1, x_2, x_3)=\phi(y_1, y_2, y_3)$} \label{Sec2}
If we take $y_1, y_2, y_3$ to be any permutation of the numbers $x_1, x_2, x_3$, we immediately get a trivial solution of Eq.~\eqref{octiceqn}. Solutions of Eq.~\eqref{octiceqn} in which  $y_1, y_2, y_3$ is not a permutation of $x_1, x_2, x_3$ will be considered nontrivial.

Since Eq.~\eqref{octiceqn} is homogeneous, any rational solution will yield, on appropriate scaling, a solution in integers. It, therefore, suffices to obtain rational  solutions of Eq.~\eqref{octiceqn}.

We will now obtain nontrivial solutions of  the diophantine equation~\eqref{octiceqn}.
We  note that if we take $x_3=0$, $y_3=0$, Eq.~\eqref{octiceqn} reduces to the equation $(x_1^2-x_2^2)^2=(y_1^2-y_2^2)^2$ for which we can readily obtain the complete solution. Such a  solution is, however, of little interest.   We also note that the form $\phi(x_1, x_2, x_3)$ is reducible over $\mathbb{Q}$, and in fact,
\begin{equation*}
\phi(x_1, x_2, x_3)=(x_1^2 - x_2^2 - x_3^2)(x_1^2 - x_2^2 + x_3^2)(x_1^2 + x_2^2 - x_3^2)(x_1^2 + x_2^2 + x_3^2).
\end{equation*}
It follows that parametric solutions of Eq.~\eqref{octiceqn} may be obtained quite simply by solving  two quadratic equations, for instance, $x_1^2 = x_2^2 + x_3^2$ and $y_1^2 = y_2^2 + y_3^2$, when both sides of \eqref{octiceqn} become $0$. Such parametric solutions of Eq.~\eqref{octiceqn} are also not very interesting since they  neither yield  the desired examples of two triangles nor do they  lead to  an elliptic curve~\eqref{cubicec} since $\phi(x_1, x_2, x_3)=0$. We therefore do not give such solutions explicitly.

To obtain parametric solutions of Eq.~\eqref{octiceqn} with $\phi(x_1, x_2, x_3) \neq 0$, we write,
\begin{equation}
\begin{aligned}
x_2&=pu, & x_3&=qv, & y_2&=pv, & y_3&=qu, \label{subsxy}
\end{aligned}
\end{equation}
where $p, q, u$ and $v$ are arbitrary parameters. On substituting these values of $x_i$, $y_i$, $i=2, 3$, in Eq.~\eqref{octiceqn}, and transposing all terms to the left-hand side, we get,
\begin{equation}
\{x_1^4 + y_1^4 - (p^4 + q^4)(u^4 + v^4)\}\{x_1^4 - y_1^4 - (p^4 - q^4)(u^4 - v^4)\}=0. \label{octiceqnred}
\end{equation}
Solutions of Eq.~\eqref{octiceqn} may be obtained by equating to $0$ either of the two factors on the left-hand side of Eq.~\eqref{octiceqnred}.

The first factor on the left-hand side of Eq.~\eqref{octiceqnred} will become $0$ if we choose $x_1, y_1, p, q, u$ and $v$ such that
\begin{equation}
(p^4 + q^4)(u^4 + v^4)=x_1^4 + y_1^4. \label{fac1}
\end{equation}
Several parametric equations of the diophantine equation~\eqref{fac1} (with a slightly different notation) have been given by Choudhry {\it et al.}~\cite{CBJ} who have also shown how infinitely many parametric solutions of Eq.~\eqref{fac1} may be obtained. These parametric solutions will yield infinitely many parametric solutions of Eq.~\eqref{octiceqn}.

As an example, a parametric solution of Eq.~\eqref{fac1},  given in \cite{CBJ}, is as follows:
\begin{equation*}
\begin{aligned}
p &= t^2 - 3,  &  q &= 4t, \\u &= (t^2 - 3)(t^2 + 9)(t^2 + 1),  &  v &= 4t(t^2 + 2t + 3)(t^2 - 2t + 3),\\
x_1 &= 16t^2(t^2 - 3)(t^2 + 3),  &  y_1 &= t^8 + 4t^6 + 86t^4 + 36t^2 + 81,
\end{aligned}
\end{equation*}
where $t$ is an arbitrary rational parameter. Using the relations~\eqref{subsxy}, we immediately get the following solution of Eq.~\eqref{octiceqn}:
\begin{equation}
\begin{aligned}
x_1 &= 16t^2(t^2 - 3)(t^2 + 3),  &   x_2 &= (t^2 - 3)^2(t^2 + 9)(t^2 + 1),\\ x_3 &= 16t^2(t^2 + 2t + 3)(t^2 - 2t + 3),  &   y_1 &= t^8 + 4t^6 + 86t^4 + 36t^2 + 81,\\ y_2 &=4t(t^2 - 3)(t^2 + 2t + 3)(t^2 - 2t + 3),  &   y_3 &= 4t(t^2 - 3)(t^2 + 9)(t^2 + 1).\\
\end{aligned}
\label{sol1octiceqn}
\end{equation}

It is readily seen that when we take $t$ to be a rational number such that $3.46 < t < 4.64$, both the sets of numbers $\{x_1^2, x_2^2, x_3^2\}$ and  $\{y_1^2, y_2^2, y_3^2\}$ satisfy the triangle inequalities, and hence we get actual examples of two triangles whose areas are equal and whose sides are perfect squares.

As a numerical example, when $t=4$, we get two triangles whose sides have lengths  $63232^2, 71825^2,  76032^2$ and $ 104593^2,  61776^2,  88400^2$, and whose areas are equal. We note that the common area of the two triangles is not rational.

We will now obtain more parametric solutions of Eq.~\eqref{octiceqn} by equating to $0$ the second factor on the left-hand side of
Eq.~\eqref{octiceqnred}. This gives a condition which, after suitable transpositions, may be written as follows:
\begin{equation}
x_1^4+hv^4 =y_1^4+hu^4, \label{fac2}
\end{equation}
where $h=p^4-q^4$. We note that a solution of Eq.~\eqref{fac2} is $(x_1, y_1, u, v)= (p, q, 1, 0)$. However,  using this solution and the relations~\eqref{subsxy}, we get only a trivial solution of Eq.~\eqref{octiceqn}.

If a solution of the diophantine equation $A^4+hB^4=C^4+hD^4$ is known such that $A \neq \pm C$ and $B \neq \pm D$, other nontrivial solutions of this equation can be found by using a method described by Choudhry~\cite{Ch}. Applying Choudhry's method to the  known solution $(x_1, y_1, u, v)= (p, q, 1, 0)$ of Eq.~\eqref{fac2}, we get the following solution of Eq.~\eqref{fac2}:
\begin{equation}
\begin{aligned}
x_1 &= p(2p^4 + 3p^3q + 3p^2q^2 + 3pq^3 - q^4),\\
y_1 &= q(p^4 - 3p^3q - 3p^2q^2 - 3pq^3 - 2q^4),\\
 u &= 2p^4 + 3p^3q + 3p^2q^2 + 3pq^3 + 2q^4, \\
v &= 3(p^2 + pq + q^2)pq,
\end{aligned}
\label{solfac2}
\end{equation}
where $p$ and $q $ are arbitrary rational parameters. We can now repeat the process taking \eqref{solfac2} as the known solution of Eq.~\eqref{fac2} and obtain a second parametric solution of Eq.~\eqref{fac2}, and by repeatedly applying this process, we can obtain a sequence of parametric solutions of Eq.~\eqref{fac2}. The parametric solution obtained by using the solution~\eqref{solfac2} is of degree $17$ in the arbitrary parameters $p$ and $q$, and subsequent solutions are of even higher degrees. Accordingly, we do not give these solutions explicitly.

The parametric solutions of Eq.~\eqref{fac2} may be used to obtain parametric solutions of the octic equation~\eqref{octiceqn}. As an example, using the solution~\eqref{solfac2} and the relations~\eqref{subsxy}, we get a solution of Eq.~\eqref{octiceqn} in terms of homogeneous polynomials of degree $5$ in the parameters $p$ and $q$. On writing $p=qt$, we may write this solution of Eq.~\eqref{octiceqn} in terms of a single rational parameter $t$ as follows:
\begin{equation}
\begin{aligned}
x_1 &= t(2t^4 + 3t^3 + 3t^2 + 3t - 1),\\
x_2 &= t(2t^4 + 3t^3 + 3t^2 + 3t + 2), \\
x_3 &= 3t(t^2 + t + 1),\\
y_1 &= t^4 - 3t^3 - 3t^2 - 3t - 2, \\
y_2 &= 3t^2(t^2 + t + 1), \\
y_3 &= 2t^4 + 3t^3 + 3t^2 + 3t + 2.
\end{aligned}
\label{sol2octiceqn}
\end{equation}

We note that when we take $t$ to be a rational number such that $ 0.53 < t < 1.88 $, both the sets of numbers $\{x_1^2, x_2^2, x_3^2\}$ and $\{y_1^2, y_2^2, y_3^2\}$ satisfy the triangle inequalities, and hence, as in the case of the parametric solution~\eqref{sol1octiceqn}, here also we get actual examples of two triangles whose areas are equal and whose sides are perfect squares.

As a numerical example, when $t=3/2$, we get, after appropriate scaling, two triangles whose sides have lengths $732^2, 804^2, 342^2$, and $293^2, 513^2, 536^2$ and whose areas are equal.

\section{Families of elliptic curves of rank~$\geq 5$} \label{Sec3}
As earlier mentioned in the Introduction, with each solution of Eq.~\eqref{octiceqn}, we may associate an elliptic curve~\eqref{cubicec} on which there are, in general, six rational points (not necessarily distinct) whose abscissae are $x_1^2x_2^2$, $x_1^2x_3^2$, $x_2^2x_3^2$, $y_1^2y_2^2$, $y_1^2y_3^2$, $y_2^2y_3^2$. It has been shown in Section~\ref{Sec2} that there exist infinitely many parametric solutions of the diophantine equation~\eqref{octiceqn}. Since \eqref{subsxy} implies $x_2x_3=y_2y_3$, each of these solutions yields a family of elliptic curves~\eqref{ecDgen} on which we know $5$ rational points.

The elliptic curve associated with the solution~\eqref{sol1octiceqn} with $t=4$ is given by
$$E:~Y^2=X^3-\frac{1617508083022593897795364438996422549375}{4}X$$
on which we get five distinct rational points given below:
\begin{align*}
P_1&=(20626479356354560000, 20849350546566884379967280000),\\
P_2&=(23113550675916619776, 54786013676180111350031745024),\\
P_3&=(29822503524802560000, -120266596559434914904320720000),\\
P_4&=(41748877998578260224, 236399461657400030514050179368),\\
P_5&=(85488908030849440000, 768253474718253155853728585000).
\end{align*}

The value of the regulator of these five points, as computed by SageMath~\cite{S}, is  $122787391.171313$. Since this is nonzero, these five points are linearly independent. It now follows from the specialization theorem of
Silverman~\cite[Theorem~11.4, p.\ 271]{S1} that the family of elliptic curves~\eqref{cubicec} related to the solutions~\eqref{sol1octiceqn} has rank at least five over $\mathbb{Q}(t)$ with linearly independent points $P_i(t)$, $i=1, \ldots, 5$.

We also get a nonzero value of the regulator for five of the rational points on the elliptic curve associated with the solution~\eqref{sol2octiceqn} with $t=3/2$. Thus, the solution~\eqref{sol2octiceqn} of Eq.~\eqref{octiceqn} also yields a family of elliptic curves of rank at least five over $\mathbb{Q}(t)$.

Computations show that the  family of elliptic curves with the solution~\eqref{sol2octiceqn} has members of higher rank in comparison to members of the  family of curves associated with the solution~\eqref{sol1octiceqn}. In the following table, we give values of $t$ which yield examples of higher rank elliptic curves for the family associated with the solution~\eqref{sol2octiceqn}:
 \begin{table}[htbp]
 \caption{\label{tab:1}High rank elliptic curves associated with the solution~\eqref{sol2octiceqn}} \vspace{-0.2cm}
 \begin{tabular}{cc}
  \toprule
  Rank          &       Values of $t$        \\
  \midrule
  $9$ &  $\displaystyle -\frac{11}{14}$,  $\displaystyle \frac{2}{9}$                   \vspace{.2cm} \\
  $8$ &  $\displaystyle -\frac{4}{21}$, $\displaystyle -\frac{12}{13}$, $\displaystyle -\frac{9}{13}$,  $\displaystyle -\frac{7}{9}$, $\displaystyle -\frac{7}{8}$, $\displaystyle -\frac{1}{8}$, $\displaystyle -\frac{3}{7}$, $\displaystyle -\frac{1}{4}$, $\displaystyle \frac{1}{12}$  \\  \bottomrule
 \end{tabular}
\end{table}

This observation motivates to study more closely  the rank of the members of the family of elliptic curves associated with the solution~\eqref{sol2octiceqn}. In the next subsection we show that the rank of the elliptic curves of this family over $\mathbb{Q}(t)$ is five and present their free generators.

\subsection{Rank and free generators of the family of elliptic curves associated with the solution~\eqref{sol2octiceqn}}
To simplify the notation, we write
$$\begin{aligned}
h_1(t) &=t^4-3t^2-2, &  h_2(t) &=2t^4+3t^2+2, \\
h_3(t) &=2t^4+3t^2-1, & h_4(t) &=4t^4+12t^3+15t^2+12t+4,   \\
h_5(t) &=4t^4+12t^3+15t^2+12t+7, & h_6(t) &=7t^4+12t^3+15t^2+12t+4.
\end{aligned}$$

By defining $(U,V)=\big((6m)^2X, (6m)^3Y\big)$ with $m=12/t^2$,
we rewrite our family of elliptic curves as
\begin{equation}
{\mathcal{E}}_t: V^2=U^3+A_4U, \label{ecfamily}
\end{equation}
where $A_4=36\prod_{i=1}^{6}h_i(t)$, with  the five
 linearly independent points being
\begin{align*}
P_1(t)&=\bigl(m^2x_1^2x_2^2, m^3x_1x_2(x_1^4+x_2^4-x_3^4)/2)\bigr),\\
P_2(t)&=\bigl(m^2x_1^2x_3^2, m^3x_1x_3(x_1^4-x_2^4+x_3^4)/2)\bigr),\\
P_3(t)&=\bigl(m^2x_2^2x_3^2, m^3x_2x_3(x_1^4-x_2^4-x_3^4)/2)\bigr),\\
P_4(t)&=\bigl(m^2y_1^2y_2^2, m^3y_1y_2(y_1^4+y_2^4-y_3^4)/2)\bigr),\\
P_5(t)&=\bigl(m^2y_1^2y_3^2, m^3y_1y_3(y_1^4-y_2^4+y_3^4)/2)\bigr).
\end{align*}

As we wish to apply a theorem of Gusi\'{c} and Tadi\'{c}~\cite[Theorem~1.3, p.\ 139]{G-T}, we write the right-hand side of \eqref{ecfamily} as follows:
$$U^3+A_4U=(U-e_1)(U-e)(U-\bar{e})$$
where
$$\begin{aligned}
e_1 &= 0, & e &= 6\sqrt{-\prod_{i=1}^{6}h_i(t)}, & \bar{e} &=-e,
\end{aligned}$$
so that
$$\begin{aligned}
e_1^2-(e+\bar{e})e_1+e\bar{e}&=36\prod_{i=1}^{6}h_i(t), & (e-\bar{e})^2&=-144\prod_{i=1}^{6}h_i(t).
\end{aligned}$$

We note that when we specialize at $t=2$, none of the nonconstant square-free divisors of $e_1^2-(e+\bar{e})e_1+e\bar{e}$ or $(e-\bar{e})^2$ in $\mathbb{Z}[t]$ is a square in $\mathbb{Q}$. Thus the family  ${\mathcal{E}}_t$ satisfies the conditions of the aforementioned theorem of Gusi\'{c} and  Tadi\'{c}, and it follows from the theorem that the specialization homomorphism 
${\mathcal{E}}_t(\mathbb{Q}(t))\rightarrow {\mathcal{E}}_{2}(\mathbb{Q})$ is injective.

Further, the value  $t=2$ leads to the rank~$5$ elliptic curve
${\mathcal E}_2$ with the torsion subgroup $\mathbb{Z}/2\mathbb{Z}$ and generators:
$$\begin{aligned}
G_1&=(680800, 1449831680), & G_2&=(981088, 1875840256), \\
G_3&=(240126016/49, 3919014764288/343), & G_4&=(55264356, -411011675928), \\
G_5&=(123121216, -1366271251712).
\end{aligned}$$

We also note that at $t=2$, the regulator of the specialized points
$$\begin{aligned}
P_1&=(123121216, 1366271251712), & P_2&=(9400356, -29246291928), \\
P_3&=(10188864, -32931327744), & P_4&=(1382976, -2504823552), \\
P_5&=(1132096, 2102770432),
\end{aligned}$$
on the specialized elliptic curve ${\mathcal{E}}_2: V^2=U^3+2624072905728U$ is the nonzero value  $123017.788734562$, showing that these points are linearly independent.
Till now, by \cite[Theorem~1.3]{G-T}, we have shown that rank~${\mathcal{E}}_t(\mathbb{Q}(t))=5$.

To determine a set of five free generators of ${\mathcal{E}}_t(\mathbb{Q}(t))$, we note the following relations between the generators and the linearly independent points:
$$\begin{aligned}
P_1&=-G_5, & P_2&=-2G_1-2G_2+G_3+G_4-G_5, \\
P_3&=-G_3-G_4, & P_4&=-2G_1-2G_2+G_4-G_5, \\
P_5&=-2G_2-G_3-G_5.
\end{aligned}$$

Since the change of basis matrix, i.e.,
$$
\begin{pmatrix}
0 & 0 & 0 & 0 & -1\\
-2 & -2 & 1 & 1 & -1\\
0 & 0 & -1 & -1 & 0\\
-2 & -2 & 0 & 1 & -1\\
0 & -2 & -1 & 0 & -1
\end{pmatrix}
$$
is not unimodular, the set $\{P_1(2), P_2(2), P_3(2), P_4(2), P_5(2)\}$ cannot be considered as a set of free generators. Indeed, the set does not generate the whole group 
${\mathcal{E}}_{2}(\mathbb{Q})/{\mathcal{E}}_{2}(\mathbb{Q})_{\rm{tors}}$, but only its subgroup of finite index. To construct a set of five free generators of ${\mathcal{E}}_t(\mathbb{Q}(t))$, we apply a method used by Dujella {\it et al.} in \cite[Section~4.2]{D-T}. Accordingly, we write
\begin{equation*}
P^{\star}_2=P_2+G_1+G_2, \quad P^{\star}_4=P_4+G_1+G_2, \quad P^{\star}_5=P_5+G_2,
\end{equation*}
and we then have
$$\begin{aligned}
P^{\star}_2&=-G_1-G_2+G_3+G_4-G_5, \\
P^{\star}_4&=-G_1-G_2+G_4-G_5, \\
P^{\star}_5&=-G_2-G_3-G_5.
\end{aligned}$$

It would be observed that the change of basis matrix is now unimodular, and hence the points $P_1$, $P^{\star}_2$, $P_3$, $P^{\star}_4$, $P^{\star}_5$ are free generators for ${\mathcal{E}}_{2}(\mathbb{Q})$. Therefore, the above injectivity is an isomorphism and hence ${\mathcal{E}}_t(\mathbb{Q}(t))$ and ${\mathcal{E}}_{2}(\mathbb{Q})$ have the same rank five, and, modulo torsion points, $P_1(t)$, $P^{\star}_2(t)$, $P_3(t)$, $P^{\star}_4(t)$, $P^{\star}_5(t)$ can be considered as free generators of ${\mathcal{E}}_t(\mathbb{Q}(t))$, where $P^{\star}_2(t)$, $P^{\star}_4(t)$, and $P^{\star}_5(t)$ are extracted as described below. For determining these points, it is sufficient to determine  the points $G_1(t)$ and $G_2(t)$. According to the relations,
$$\begin{aligned}
2G_1(t) &= -P_3(t)-P_4(t)+P_5(t), &
2G_2(t) &= P_1(t)-P_2(t)+P_4(t)-P_5(t),
\end{aligned}$$
and using the double point formula, we can get, for example,
$$\begin{aligned}
G_1(t)&=\left(\frac{2h_2(t)h_6(t){\theta}_1^2(t)}{t^4},
\frac{4h_2(t)h_6(t){\theta}_1(t){\theta}_2(t)}{t^{6}}\right),\\
G_2(t)&=\left(\frac{4h_1(t)h_2(t)h_3(t)h_4(t)}{t^2}, \frac{4h_1(t)h_2(t)h_3(t)h_4(t){\theta}_3(t)}{t^3}\right),
\end{aligned}$$
where
$$\begin{aligned}
{\theta}_1(t)&=2t^4-3t^3-3t^2-6t-4, \\
{\theta}_2(t)&=16t^{12}+12t^{11}+66t^{10}+27t^9-27t^8-183t^7-381t^6-522t^5-588t^4\\
&\quad -492t^3-312t^2-144t-32,\\
{\theta}_3(t)&=8t^8+12t^7+6t^6+3t^5-20t^4+3t^3+6t^2+12t+8.
\end{aligned}$$

With these values of $G_1(t)$ and $G_2(t)$,  a set of free generators of ${\mathcal{E}}_t(\mathbb{Q}(t))$ is given by $P_1(t)$, $P^{\star}_2(t)$, $P_3(t)$, $P^{\star}_4(t)$, $P^{\star}_5(t)$. We note that this set of generators is not necessarily unique.

We have  verified that a similar  result holds for the first family. Since the proof is similar, we omit the details.

\section{Concluding remarks} \label{Sec4}
In this work, we studied the diophantine equation $\phi(x_1, x_2, x_3)=\phi(y_1, y_2, y_3)$ with
$\phi(x_1, x_2, x_3)=x_1^8+ x_2^8 + x_3^8 - 2x_1^4x_2^4 - 2x_1^4x_3^4- 2x_2^4x_3^4$, which from the geometric point of view is equivalent to two equiareal triangles with squared sides, and to each solution we linked a family of elliptic curves of rank~$5$.  In 2021, the authors~\cite{Ch-Sh1} obtained parametric as well as numerical solutions of the sextic diophantine chain $\varphi(x_1, x_2, x_3) = \varphi(y_1, y_2, y_3) = \varphi(z_1, z_2, z_3) $
with $\varphi(x_1, x_2, x_3) = x_1^6 + x_2^6 + x_3^6 - 2x_1^3x_2^3 - 2x_1^3x_3^3 - 2x_2^3x_3^3$ and
constructed a parametrized family of Mordell curves $y^2 = x^3 + \varphi(x_1, x_2, x_3)/4$ of generic rank~$\geq 6$ using the parametric solution of the sextic diophantine chain. So, one can naturally pose the problem of
determining other positive integer values of $k$ and $\ell$ for which
the system of equations
$$\phi_{\ell}(x_{11}, x_{12}, x_{13})=\phi_{\ell}(x_{21},x_{22}, x_{23})=\cdots=\phi_{\ell}(x_{k1}, x_{k2}, x_{k3})$$
has rational solutions, where
$$\phi_{\ell}(x_{11}, x_{12}, x_{13})=
x_{11}^{2\ell}+x_{12}^{2\ell}+x_{13}^{2\ell}-2x_{11}^{\ell}x_{12}^{\ell}-2x_{11}^{\ell}x_{13}^{\ell}-2x_{12}^{\ell}x_{13}^{\ell}.$$

It  would  also be of interest to find more than two equiareal triangles with squared sides.


\begin{thebibliography}{10}

\bibitem{Ch} A. Choudhry, \textit{On the diophantine equation $A^4+hB^4=C^4+hD^4$}, Indian J.  Pure Appl. Math. {\bf  26}  (1995), 1057--1061.

\bibitem{CBJ} A. Choudhry, I. Bluskov and A. James, \textit{A diophantine equation inspired by Brahmagupta's identity}, Int. J. Number Theory {\bf 18} (2022), 905--911.

\bibitem{Ch-Sh1} A. Choudhry and A. Shamsi Zargar, \textit{A sextic diophantine chain and a related Mordell curve}, Integers {\bf 21} (2021), $\#$ 39.

\bibitem{Ch-Sh} A. Choudhry and A. Shamsi Zargar, \textit{Pairs of equiperimeter and equiareal triangles whose sides are perfect squares}, Func. Funct. Approx. Comment. Math. (to appear).

\bibitem{D-T} A. Dujella, J. C. Peral and P. Tadi\'{c}, \textit{Elliptic curves with torsion group $\mathbb{Z}/6\mathbb{Z}$}, Glas. Mat. Ser. III {\bf 51} (2016), 321--333.

\bibitem{G-T} I. Gusi\'{c}, P. Tadi\'{c}, \textit{Injectivity of specialization homomorphism of elliptic curves}, J. Number Theory {\bf 148} (2015), 137--152.

\bibitem{IN} F. Izadi and K. Nabardi, \textit{A family of elliptic curves with rank~$\geq5$}, Period. Math. Hungar. {\bf 71} (2015), 243--249.

\bibitem{Ma} M. Maenishi, \textit{On the rank of elliptic curves $y^2=x^3-pqx$}, Kumamoto J. Math. {\bf 15} (2002), 1--5.

\bibitem{Nag} K.-I Nagao, \textit{On the rank of elliptic curve $y^2=x^3-kx$}, Kobe J. Math. {\bf 11} (1994), 205--210.

\bibitem{S} Sage Developers, \textit{SageMath, Sage mathematics software system}, available at \url{https://www.sagemath.org}.

\bibitem{S1} J. H. Silverman, \textit{Advanced Topics in the Arithmetic of Elliptic Curves},
Springer, New York, 1994.

\end{thebibliography}
\end{document}